\documentclass[12pt]{amsart}

\usepackage{amsmath}
\usepackage{amscd}
\usepackage{amssymb}
\usepackage{latexsym}
\usepackage{stmaryrd}
\usepackage{mathbbol}
\usepackage{mathabx}
\usepackage{color}

\input xy
\xyoption{all} \CompileMatrices

\definecolor{orange}{rgb}{1,0,0}
\definecolor{cadmiumgreen}{rgb}{0, 0.78, 0.05}

\newtheorem{theorem}{Theorem}[section]
\newtheorem*{theorem*}{Theorem}
\newtheorem{lemma}[theorem]{Lemma}
\newtheorem*{lemma*}{Lemma}
\newtheorem{corollary}[theorem]{Corollary}
\newtheorem{proposition}[theorem]{Proposition}
\newtheorem{defprop}[theorem]{Proposition-Definition}
\newtheorem{remark}[theorem]{Remark}
\newtheorem{definition}[theorem]{Definition}

\newtheorem{propdef}[theorem]{Proposition-Definition}
\makeatletter
\def\revddots{\mathinner{\mkern1mu\raise\p@
\vbox{\kern7\p@\hbox{.}}\mkern2mu
\raise4\p@\hbox{.}\mkern2mu\raise7\p@\hbox{.}\mkern1mu}}
\makeatother 
\newcommand{\bgl}{\begin{equation}} 
\newcommand{\egl}{\end{equation}}
\newcommand{\bgloz}{\begin{equation*}} 
\newcommand{\egloz}{\end{equation*}}
\newcommand{\bgln}{\begin{eqnarray}} 
\newcommand{\egln}{\end{eqnarray}}
\newcommand{\bglnoz}{\begin{eqnarray*}} 
\newcommand{\eglnoz}{\end{eqnarray*}}
\newcommand{\btheo}{\begin{theorem}}
\newcommand{\etheo}{\end{theorem}}
\newcommand{\btheooz}{\begin{theorem*}}
\newcommand{\etheooz}{\end{theorem*}}
\newcommand{\blemma}{\begin{lemma}}
\newcommand{\elemma}{\end{lemma}}
\newcommand{\blemmaoz}{\begin{lemma*}}
\newcommand{\elemmaoz}{\end{lemma*}}
\newcommand{\bproof}{\begin{proof}}
\newcommand{\eproof}{\end{proof}}
\newcommand{\bbew}{\begin{beweis}}
\newcommand{\ebew}{\end{beweis}}
\newcommand{\bremark}{\begin{remark}\em}
\newcommand{\eremark}{\end{remark}}
\newcommand{\bdefin}{\begin{definition}}
\newcommand{\edefin}{\end{definition}}
\newcommand{\bprop}{\begin{proposition}}
\newcommand{\eprop}{\end{proposition}}
\newcommand{\bdefprop}{\begin{defprop}}
\newcommand{\edefprop}{\end{defprop}}
\newcommand{\bcor}{\begin{corollary}}
\newcommand{\ecor}{\end{corollary}}
\newcommand{\bfa}{\begin{cases}} 
\newcommand{\efa}{\end{cases}}
\newcommand{\nn}{\par\vspace{-3mm}\noindent}

\newcommand{\cK}{\mathcal K}

\newcommand{\cM}{\mathcal M}

\def\Nz{\mathbb{N}}

\def\Zz{\mathbb{Z}}

\def\1z{\mathbb{1}}
\newcommand{\vp}{\varphi}

\def\SEMI{\mbox{$\times\kern-2pt\vrule height5pt width.6pt \kern3pt $}}

\newcommand{\id}{{\rm id}}

\newcommand{\Ad}{{\rm Ad\,}}

\newcommand{\ev}{\operatorname{ev}} 

\begin{document}
\title[Quasihomomorphisms, operator homotopy and stable uniqueness]{Split exactness, operator homotopy and stable uniqueness in $KK$}

\author{Joachim Cuntz}
\address{Mathematisches Institut\\
 Universit\"at M\"unster\\
  Ein\-stein\-str.\ 62\\
  48149 M\"unster\\
  Germany}
\email{cuntz@uni-muenster.de}


\begin{abstract}
  We develop important properties of the $KK$-functor on the basis of split exactness.
\end{abstract}
\thanks{Supported by Deutsche Forschungsgemeinschaft (DFG) via Exzellenzstrategie des Bundes und der L\"ander EXC 2044–390685587, Mathematik M\"unster: Dynamik–Geometrie–Struktur.}
\maketitle
\section{Introduction}
The bivariant functor $KK$ introduced by Kasparov is a fundamental tool in the theory of C*-algebras.
In our recent paper \cite{CuGab} it was remarked that there is a simple approach to the product in $KK$ on the basis of the split exactness property. On the other hand it is known since a long time that homotopy invariance is a consequence of split exactness together with stability \cite{HigAlgK}. There is a more recent rather simple proof of this fact by Ralf Meyer \cite[Lemma 3.26]{CMR}.\\
Thus key properties of $KK$ are direct consequences of split exactness. In the present note we explain and develop these two points on the basis of the framework of quasihomomorphisms introduced in \cite{CuGen},\cite{CuKK}.
In Section 5 we first construct the product using the framework of the universal algebra $qA$. An important ingredient is a quasihomomorphism associated with a split exact sequence $0\to J\to A\to B\to 0$ of C*-algebras which represents the $KK$-version of a projection onto the ideal in this sequence. This leads in particular to a somewhat more direct construction of the universal map $\vp_A:qA\to M_2(q^2A)$ of \cite{CuKK,CuGab} and to a somewhat simpler proof of the properties of $\vp_A$ that are used in \cite{CuKK,CuGab} to prove the associativity of the product $KK(A,B)\times KK(B,C)\to KK(A,C)$.\\
We then turn in Section 6 to a slightly different description of $KK$ in terms of quasihomomorphisms (without using $qA$) and also discuss operator homotopy. We establish the product and its associativity for the version $KK^{op}$ of $KK$ which is defined using operator homotopy. We then give a short proof for the fact that $KK^{op}$ is homotopy invariant and thus coincides with $KK$.
It turns out that Meyer's proof for homotopy invariance becomes even shorter in our setting. We thus obtain a simple proof of Kasparov's homotopy invariance result that shows that operator homotopy together with unitary equivalence of Kasparov modules gives homotopy. As a consequence, we show in Corollary \ref{cordeg} that two quasihomomorphisms $(\alpha,\bar{\alpha})$
and $(\beta,\bar{\beta})$ from $A$ to $\cK\otimes B$ define the same element in $KK(A,B)$ iff
we can add degenerates to $(\alpha,\bar{\alpha})$
and $(\beta,\bar{\beta})$ such that the sums become operator homotopic. This corollary is important
for applications.\\
Our treatment of $KK$ in sections 6-8 and the relatively short proof of \ref{cordeg} is independent of Section 5 that uses the $qA$-formalism. Thus we get two closely related and similar, but different, approaches to $KK$ - the second approach being more suitable for the discussion of $KK^{op}$.\\
A version of Corollary \ref{cordeg} that is especially useful and important in the classification program for nuclear C*-algebras is the so called stable uniqueness theorem \cite{DadEil}, \cite{Lin}. For its applications see e.g. \cite{DadEilC}, \cite{Lin},\cite{Gabe},\cite{Schafhauser},\cite{TWW}. It shows that operator homotopy can be replaced by an even more restrictive equivalence relation and that a quasihomomorphism $(\alpha,\bar{\alpha})$
from $A$ to $\cK\otimes B$ represents the element 0 in $KK(A,B)$ iff we can add a degenerate to $(\alpha,\bar{\alpha})$ such that $\alpha$ and $\bar{\alpha}$ become asymptotically unitarily homotopic with unitaries $U_t$ in $1+ \cK\otimes B$. Already in \cite{DadEil} this theorem was deduced from the result on operator homotopy. But there is a more direct way that leads from the operator homotopy result to the stable uniqueness theorem using quasicentral approximate units. I am grateful to Gabor Szabo who pointed out the argument in \cite[Lemma 4.3]{GabSzab}. In section \ref{secstab} we use the main idea in \cite[Lemma 4.3]{GabSzab} to give a very short proof of the stable uniqueness theorem (in the non-equivariant case). Combining this with our short proof of Corollary \ref{cordeg} we obtain a self-contained simple proof of the stable uniqueness theorem. \\
The idea to use the split exactness property for the discussion of the Kasparov product is of course not new, cf. e.g. \cite{HigChar},\cite[17.8.4]{BB},\cite[18.11.1]{BB} and goes back to \cite{CuGen}. But here we make use of this property more systematically. \\
Using the approach in \cite{CuGab} the arguments in this paper could also be extended to some of the versions of $KK$ with extra structures studied in \cite{CuGab} but we don't do that here. In particular, such an extension to $KK^{nuc}$ in place of $KK$ should be straightforward. I am also indebted to Marius Dadarlat for helpful remarks on an earlier version of this manuscript.
\section{Preliminaries}\label{s0}
Notation: In the following, homomorphisms between C*-algebras will always be assumed to be  *-homomorphisms. By $\cK$ we denote the standard algebra of compact operators on $\ell^2 \Nz$. There is a natural isomorphism $\cK\cong \cK\otimes \cK$. A C*-algebra $A$ is called stable if $A\cong \cK\otimes A$. Given a C*-algebra $A$ we denote by $\cM (A)$ its multiplier algebra. If $\vp:A\to B$ is a $\sigma$-unital homomorphism between C*-algebras, we denote by $\vp^\circ$ its extension to a homomorphism $\cM (A)\to \cM (B)$.\\
Let $A$ be a C*-algebra. We denote by $QA$ the free product $A\star A$ and by $\iota, \bar{\iota}$ the two natural inclusions of $A$ into $QA=A\star A$. We denote by $qA$ the kernel of the natural map $A\star A\to A$ that identifies the two copies $\iota (A)$ and $\bar{\iota} (A)$ of $A$. Then $qA$ is the closed two-sided ideal in $QA$ that is generated by the elements $qx=\iota (x) - \bar{\iota} (x),\,x\in A$. The characteristic identity for the map $x\mapsto qx$ is the identity $q(xy)=qxqy - \iota (x)qy - qx  \,\iota (y)$. More generally, if $\vp$ is a homomorphism and $\delta$ a linear map between two algebras, then $\vp -\delta$ is an algebra homomorphism if and only if $\delta$ satisfies
\bgl\label{formq} \delta(xy)=\delta x\delta y - \vp(x)\delta y - \delta x  \,\vp (y) \egl \nn
There is the natural evaluation map $\pi_A: qA \to A$ given by the restriction to $qA$ of the map $\id\star 0:QA\to A$ that is the identity on the first copy of $A$ and zero on the second one. Similarly there is a second evaluation map $\breve{\pi}_A$ which evaluates at the second copy.
As in \cite{CuKK} we define a prequasihomomorphism between two C*-algebras $A$ and $B$ to be a diagram of the form \bgl\label{fpq}A \quad \stackrel{\vp,\bar{\vp}}{\rightrightarrows}\quad E \rhd J \stackrel{\mu}{\to} B\egl
i.e. two homomorphisms $\vp,\bar{\vp}$ from $A$ to a C*-algebra $E$ that contains an ideal $J$, with the condition that $\vp (x)-\bar{\vp} (x)\in J$ for all $x\in A$ and finally a homomorphism $\mu : J\to B$. The pair $(\vp,\bar{\vp})$ induces a homomorphism $QA\to E$ by mapping the two copies of $A$ via $\vp,\bar{\vp}$. This homomorphism maps the ideal $qA$ to the ideal $J$. Thus, after composing with $\mu$, every such prequasihomomorphism from $A$ to $B$ induces naturally a homomorphism $q(\vp,\bar{\vp}):qA\to B$. Conversely, if $\psi: qA \to B$ is a homomorphism, then we get a prequasihomomorphism by choosing $E =\cM(\psi (qA)), \,J=\psi(qA)$ and $\vp=\psi^\circ \iota,\,\bar{\vp}=\psi^\circ \bar{\iota}$ as well as the inclusion $\mu: \psi(qA) \hookrightarrow B$.\\
Moreover, $A\mapsto qA$ is a functor. Any homomorphism $\alpha:A\to B$ induces a homomorphism $q \alpha:qA\to qB$ that maps $q(x)$ to $q(\alpha x)$. In the notation above it can be described as $q(\iota_B\alpha,\bar{\iota}_B\alpha)$.
We say that a prequasihomomorphism is a quasihomomorphism if the map $\mu$ is simply an inclusion (this convention differs slightly from the notation in \cite{CuKK}).\\ If a quasihomomorphism is represented by the diagram \eqref{fpq} with $\mu$ an inclusion, we will use the shorthand notation $(\vp,\bar{\vp})$ for this quasihomomorphism.

\section{Quasihomomorphisms and $KK$}
Given two homomorphisms $\vp,\psi:X\to Y$ between C*-algebras we denote by $\vp\oplus\psi$ the homomorphism $$x\mapsto\;\scriptsize{\left(\begin{matrix}
 \vp (x) & 0\\
 0 & \psi (x)
\end{matrix}\right)}$$ from $X$ to $M_2(Y)$. Following \cite{CuKK} we define
\bdefin
Let $A$, $B$ be C*-algebras and $qA$ as in Section \ref{s0}. We define $KK(A,B)$ as the set of homotopy classes of homomorphisms from $qA$ to $\cK\otimes B$.
\edefin
The set $KK(A,B)$ becomes an abelian group with the operation $\oplus$ that assigns to two homotopy classes $[\alpha],[\beta]$ of homomorphisms $qA\to \cK\otimes B$ the homotopy class $[\alpha\oplus \beta]$ (using an isomorphism $M_2(\cK) \cong \cK$ to identify $M_2(\cK\otimes B)\cong \cK\otimes B$; this is well-defined since such an isomorphism is unique up to homotopy). The additive inverse to $\alpha =q(\vp,\bar{\vp})$ is $\breve{\alpha} =q(\bar{\vp},\vp)$. In \cite{CuGen} it was checked that this definition of $KK(A,B)$ is equivalent to the one by Kasparov.\\
A fundamental example of a $KK$-element arises as follows: Let $0\to J\to E\stackrel{p}{\to} A\to 0$ be an exact sequence of C*-algebras with a splitting $s:A\to E$. Then the quasihomomorphism $\kappa_E=(\id_E,sp)$ defines an element of $KK(E,J)$ which plays the role of a projection $E\to J$.
\bremark\label{remstab}
Every homomorphism $q(\vp,\bar{\vp}):qA\to \cK\otimes B$ can be extended naturally to a homomorphism $q(\id_\cK\otimes \vp,\id_\cK\otimes \bar{\vp}):q(\cK\otimes A)\to \cK\otimes\cK\otimes B\cong \cK\otimes B$.
\eremark
\section{Operator homotopy and $KK^{op}$}\label{secophom}
In this section we consider only quasihomomorphisms into algebras that are stabilized by $\cK$.
\bdefin
We say that a quasihomomomorphism $(\vp,\bar{\vp})$ from $A$ to $\cK\otimes B$ is $\cK$-stable if it is of the form $A \, \stackrel{\vp,\bar{\vp}}{\rightrightarrows}\, E \,\rhd J \stackrel{\mu}{\to} \cK\otimes B$ with $J=\cK\otimes B,\, E=\cM(\cK\otimes B)$ and $\mu=\id_{\cK\otimes B}$.
\edefin

\bdefin\label{defop}
Let $(\vp_1,\bar{\vp}_1)$,$(\vp_2,\bar{\vp}_2)$ be two $\cK$-stable quasihomomorphisms from $A$ to $\cK \otimes B$. \\
(a) We say that $(\vp_1,\bar{\vp}_1)$ and $(\vp_2,\bar{\vp}_2)$ are elementary operator homotopic if $\vp_1=\vp_2$ and there is a norm continuous family $[0,1]\ni t\mapsto U_t$ of unitaries in $\cM(\cK\otimes B)$ such that $U_t\bar{\vp}_1(x)-\bar{\vp}_2(x)U_t\in \cK\otimes B$ for all $x,t$ and such that, with $U=U_1$, we have $\bar{\vp}_2=\Ad U \bar{\vp}_1$.\\
(b) We say that $(\vp_1,\bar{\vp}_1)$ and $(\vp_2,\bar{\vp}_2)$ are unitarily equivalent if there is a unitary $U$ in $\cM(\cK\otimes B)$ such that $(\vp_2,\bar{\vp}_2)=(\Ad U \vp_1,\Ad U \bar{\vp}_1)$.
Replacing $U$ by $U\oplus U^*$ in $M_2(\cM(\cK \otimes B))$ we may assume that $U$ is homotopic to $1$.\\
(c) We say that $(\vp_1,\bar{\vp}_1)$ and $(\vp_2,\bar{\vp}_2)$ are operator homotopic if they can be joined by a sequence of elementary operator homotopies and unitary equivalences.
\edefin
\blemma\label{lemuni}
Two $\cK$-stable quasihomomorphisms $(\vp_1,\bar{\vp}_1)$ and $(\vp_2,\bar{\vp}_2)$ from $A$ to $\cK \otimes B$ are operator homotopic if and only if there are unitaries $V,W$ in $\cM(\cK\otimes B)$ such that $(\vp_2,\bar{\vp}_2) = (\Ad V\vp_1,\Ad W\bar{\vp}_1)$ and there is a norm continuous family $[0,1]\ni t\mapsto U_t$ of unitaries in $\cM(\cK\otimes B)$ such that $U_t\bar{\vp}_2(x)-\bar{\vp}_2(x)U_t\in \cK\otimes B$ for all $x,t$ and such that $U_0=1$ and $U_1=VW^*$.\\
In particular any operator homotopy is the composition of one unitary equivalence (induced by $\Ad V$) and one elementary operator homotopy (given by the family $U_t$).\\
Consequently $(\vp,\bar{\vp})$ is operator homotopic to $(\vp,\vp)$, iff $\vp$ is elementary operator homotopic to $\bar{\vp}$.
\elemma
\bproof
It is clear that $(\vp_2,\bar{\vp}_2) = (\Ad V\vp_1,\Ad W\bar{\vp}_1)$ arises from $(\vp_1,\bar{\vp}_1)$ by the composition of the unitary equivalence induced by $\Ad V$ and the elementary operator homotopy induced by $U_t$. Conversely it is also clear that, if we apply a unitary equivalence or an elementary operator homotopy to $(\vp_2,\bar{\vp}_2) = (\Ad V\vp_1,\Ad W\bar{\vp}_1)$, then we get a quasihomomorphism of the same form $(\Ad V'\vp_1,\Ad W'\bar{\vp}_1)$ with $V',W'$ such that $V'W'^*$ is homotopic to $1$ via a homotopy $U_t$ that commutes with $\bar{\vp}_2(x)$ mod $\cK\otimes B$ for all $x$.
\eproof
\bdefin\label{defoh}
Let $A$ and $B$ be C*-algebras. Given a $\cK$-stable quasihomomorphism $(\vp,\bar{\vp})$ from $A$ to $\cK\otimes B$ we denote by $[(\vp,\bar{\vp})]_{oph}$ the equivalence class of $(\vp,\bar{\vp})$ for the relation of operator homotopy.   We denote by $S^{op}(A,B)$ the set of equivalence classes of quasihomomorphisms $A\to \cK\otimes B$. This is an abelian semigroup with addition induced by $\oplus$. We denote by $S_0^{op}(A,B)$ the subsemigroup of equivalence classes represented by degenerate elements of the form $(\theta,\theta)$. We define $KK^{op}(A,B)$ as the quotient $S^{op}(A,B)/S_0^{op}(A,B)$.
\edefin
By definition of the quotient of an abelian semigroup $S$ by a subsemigroup $S_0$, two elements $x,y\in S$ become equal in $S/S_0$ iff there are $a,b\in S_0$ such that $x+a=y+b$.\\ Given  a quasihomomorphism $(\vp,\bar{\vp})$ from $A$ to $\cK\otimes B$ the quasihomomorphism
$$(\vp,\bar{\vp})\oplus (\bar{\vp},\vp) \,=\,
\left(\scriptsize{\left(\begin{matrix}
 \vp & 0\\
 0 & \bar{\vp}
\end{matrix}\right)}
,\,
\scriptsize{\left(\begin{matrix}
 \bar{\vp}  & 0\\
 0 & \vp
\end{matrix}\right)}\right)
$$
is elementary operator homotopic, via a rotation of the second matrix exchanging the lower right with the upper left position, to the degenerate quasihomomorphism $((\vp\oplus \bar{\vp}), (\vp\oplus \bar{\vp}))$. Therefore, if we denote by $[(\vp,\bar{\vp})]_{op}$ the image of $[(\vp,\bar{\vp})]_{oph}$ in $KK^{op}(A,B)$ then $KK^{op}(A,B)$ becomes an abelian group with the addition
$$[(\vp,\bar{\vp})]_{op} + [(\psi,\bar{\psi})]_{op} = [(\vp,\bar{\vp})\oplus (\psi,\bar{\psi})]_{op}$$
Two $\cK$-stable quasihomomorphisms $(\vp,\bar{\vp})$ and $(\psi,\bar{\psi})$ represent the same element of $KK^{op}$ i.e. $[(\vp,\bar{\vp})]_{op}=[(\psi,\bar{\psi})]_{op}$ if and only if there are degenerates $(\theta_1,\theta_1)$ and $(\theta_2,\theta_2)$ such that $[(\vp,\bar{\vp})\oplus (\theta_1,\theta_1)]_{oph}=[(\psi,\bar{\psi})\oplus (\theta_2,\theta_2)]_{oph}$. We can add a degenerate $(\theta,\theta)$ to both sides where $\theta$ contains $\theta_1$ and $\theta_2$ with infinite multiplicity and thus assume that $\theta_1=\theta_2$. It is clear that then the two quasihomomorphisms represent homotopic homomorphisms $qA\to \cK\otimes B$ and thus the same element in $KK(A,B)$.
\bremark\label{propdeg}
One could also define $KK^{op}(A,B)$ as the set of equivalence classes of homomorphisms $qA\to \cK\otimes B$ where two such homomorphisms $\vp$ and $\psi$ are called equivalent if they can be represented by two $\cK$-stable quasihomomorphisms $(\alpha,\bar{\alpha})$ and $(\beta,\bar{\beta})$ (i.e. $\vp=q(\alpha,\bar{\alpha})$, $\psi=q(\beta,\bar{\beta})$) which are operator homotopic up to addition of degenerates. In Section \ref{secprod} we will discuss the product for $KK$ on the basis of its definition via homomorphisms from $qA$ to $\cK\otimes B$. The arguments in this section would then be enough to show that the alternative definition of $KK^{op}$ also leads to homotopy invariance and to show that the natural map from the alternative $KK^{op}(A,B)$ to $KK(A,B)$ is an isomorphism. It is however not clear if in this way we could get Corollary \ref{cordeg} which is important for applications.
\eremark

\section{Split exact sequences and the product}\label{secprod}
Because we will use Thomsen's extension theorem (a variant of Kasparov's technical theorem), we will assume in this section and later that the algebras $A$ and $I$ are separable.
\bprop\label{lemext}(cf. \cite[18.11.1]{BB}, where the following statement is proved using the Kasparov product - by contrast here we use it to establish the product).
Let $0 \to I\stackrel{j}{\to} E \stackrel{p}{\to} A \to 0$ be an exact sequence of C*-algebras with a splitting homomorphism $s:A\to E$ (i.e. $ps=\id_A$) and $\beta :qI\to B$ a homomorphism. There is a homomorphism $\beta'=:qE \to M_2(B)$ so that the restriction $\beta' q(j)$ of $\beta'$ to $qI$ is homotopic to $\beta\oplus 0$.
\eprop
\bproof
Let $B_0$ denote the hereditary subalgebra of $B$ generated by the image $\beta (qI)$ and let $\beta^\circ:\cM(qI)\to \cM(B_0)$ be the extension of $\beta$ to multipliers. Composing the natural map from $E$ to the multipliers of $I$ with $\iota_I^\circ,\bar{\iota}_I^\circ$ we get two natural maps $\iota_E,\bar{\iota}_E:E\to \cM(B_0)$  such that $\iota_E j=\iota_I$ and $\bar{\iota}_E j=\bar{\iota}_I$. Let $\eta =\beta^\circ \iota_E, \bar{\eta}=\beta^\circ \bar{\iota}_E\,:E\to \cM(B_0)$. We get $\beta =q(\eta j,\bar{\eta}j)= q(\eta,\bar{\eta})q(j)$.\\
Consider the C*-algebra $R$ generated in $M_2(\cM (B_0))$ by the matrices $$\left(\begin{matrix}
 R_1 & R_1R_2\\
 R_2R_1 & R_2
\end{matrix}\right)$$
where $R_1=\eta(I)$, $R_2=\bar{\eta}(I)$. \\
Consider also the C*-subalgebra $D$ of $M_2(\cM (B_0))$ consisting of matrices of the form
\bgl\label{D} d(x)=\left(\begin{matrix}
 \eta (x) & 0\\
 0 & \bar{\eta}(x)
\end{matrix}\right)\quad x\in E\egl
Then $R$ is a subalgebra of $M_2(\cM(B_0))$. Let $J= R\cap M_2(B_0)$. Since $B_0$ is an ideal in $\cM(B_0)$ this is an ideal in $R$. One also clearly has $DR, RD\subset R$. Thus $R$ is an ideal in $R+\!D$ and $J$ is also an ideal of $R+\!D$. The quotient $(R+D)/J$ is isomorphic to $M_2(\dot{\eta} (I))+\dot{D}$ where $\dot{\eta} (I)$ is the image of $\eta(I)$ in the quotient by $J$ and $\dot{D}$ is the quotient $D/J$ and consists of the matrices in $\cM(J)/J$ of the form $\left(\scriptsize{\begin{matrix}
 \dot{\eta} (x) & 0\\
 0 & \dot{\eta} (x)
\end{matrix}}\right)$ with $x\in E$ (note that $\eta (x)$ and $\bar{\eta} (x)$ become equal as multipliers of $\dot{\eta} (I)$ for $x\in E$).\\
We can now apply Thomsen's extension theorem \cite[1.1.26]{JeTh} and lift the multiplier $\left(\scriptsize{\begin{matrix}
 0 & 1\\
 1 & 0
\end{matrix}}\right)$ of $(R+D)/J$ to a selfadjoint multiplier $S$ of $R+D$, and thus of $J$, that commutes mod $J$ with $D$. We set $F=e^{\frac{\pi i}{2} S}$ and denote by $\sigma$ the automorphism $\Ad F$ of $\cM(J)$. \\
Consider the pair of homomorphisms $E\to \cM(J)$
\bgl\label{form'}
\eta'=\left(\begin{matrix}
 \eta &0\\
 0 & \bar{\eta}\circ(sp)
\end{matrix}\right),\, \bar{\eta}'=\sigma\,\left(\begin{matrix}
 \eta\circ(sp)&0\\
 0 & \bar{\eta}
\end{matrix}\right)\egl
Note that $\eta'(x)=d(x)+(0\oplus\bar{\eta}(sp(x)-x)$. Using the fact that $\sigma$ by definition fixes $d(x)$ mod $J$ and moves $(0\oplus \bar{\eta}(sp(x)-x))$ to $\eta(sp(x)-x) \oplus 0$ mod $J$ we see that $\eta'(x) = \bar{\eta}'(x)$ mod $J$ for $x\in E$.
Therefore the pair $(\eta',\bar{\eta}')$ defines a map $\beta': qE\to J$.
Since $sp$ is 0 on $I$, the restriction of $\beta'$ to $qI$ is given by the pair
$$
\left(\left(\begin{matrix}
 \eta &0\\
 0 & 0
\end{matrix}\right)\,,\, \sigma\,\left(\begin{matrix}
 0&0\\
 0 & \bar{\eta}
\end{matrix}\right)\right)$$
and, as a homomorphism from $qI$ to $J$, is clearly homotopic to $q(\eta, \bar{\eta})\oplus 0$
\eproof
\bremark
The reader may have noticed that $(\beta',\bar{\beta}')$ is nothing but the product of $(\beta,\bar{\beta})$ with the natural projection quasihomomorphism $(\id_E,sp):E\to I$.
\eremark
We now want to define the product of $KK$-elements given by  $\alpha :qA\to \cK\otimes B$ and $\beta :qB\to \cK\otimes C$. Using Remark \ref{remstab} $\beta$ extends to a stabilized homomorphism still denoted by $\beta$ from $q(\cK\otimes B)$ to $\cK\otimes C$. For the product only the restriction of $\beta$ to $qB_0$ will matter, where $B_0$ is as above.

We define $\alpha_E,\bar{\alpha}_E:A \to \cM(B_0)\oplus A$ by $\alpha_E(x)= (\alpha^\circ\iota_A (x)\,,\, x)$, $\bar{\alpha}_E(x)= (\alpha^\circ\bar{\iota}_A (x)\,,\, x)$ and set $E_\alpha = C^*((B_0\oplus 0),\alpha_E(A))$. This gives an exact sequence
$0 \to B_0\to E_\alpha \stackrel{p}{\to} A \to 0$ with two splittings given by $\alpha_E,\bar{\alpha}_E : A\to E_\alpha$. Note that the quasihomomorphism $(\alpha_E,\bar{\alpha}_E)$ represents $\alpha:qA\to B_0$ i.e. $\alpha\oplus 0 = q(\alpha_E,\bar{\alpha}_E)$.
Recall that for a homomorphism $\mu:qA\to \cK\otimes B$ given by the pair $(\vp,\bar{\vp})$ the homomorphism $\breve{\mu} =q(\bar{\vp},\vp)$ is an additive homotopy inverse.
If $\nu$ is a second such inverse to $\mu$, then $\nu$ is homotopic to $\breve{\mu}$ in matrices (because $\nu \sim \nu\oplus\mu\oplus \breve{\mu}\sim 0\oplus 0\oplus \breve{\mu}$).
\bprop\label{propomega}
Let $\alpha$ and $\beta =q(\eta,\bar{\eta})$ as well as $B_0,E_\alpha$ be as above and assume that $\beta'=q(\eta',\bar{\eta}'):qE_\alpha \to \cK\otimes C$ extends $\beta$ up to homotopy as in \ref{lemext}.
If $C_0$ denotes the hereditary subalgebra generated by $\beta'(qE_\alpha)$ in $\cK\otimes C$, $\eta',\bar{\eta}'$ are homomorphisms $E_\alpha \to \cM(C_0)$ which we can compose with $\alpha_E, \bar{\alpha}_E:A\to E_\alpha$.\\
The homomorphism $\beta q(\alpha):q^2A\to C_0\subset \cK\otimes C$ is homotopic to $\omega q(\pi_A)$ where $\omega :qA\to C_0\subset\cK\otimes C$ is given by $\omega = \beta' q(\alpha_E)\oplus\breve{\beta}' q(\bar{\alpha}_E)$. The map $\omega$ can, more explicitly, also be described by the formula $\omega =q(\eta'\alpha_E\oplus\bar{\eta}'\bar{\alpha}_E \, ,\, \bar{\eta}'\alpha_E\oplus\eta'\bar{\alpha}_E )$.\\
\eprop
\bproof
The homomorphism $\alpha = q(\alpha_E,\bar{\alpha}_E):qA\to B_0$ extends to the homomorphism $\alpha_E\star\bar{\alpha}_E$ from $QA$ to $E_\alpha$. As a homomorphism to $M_2(E_\alpha)$ this extended map is homotopic, via a rotation in $M_2(E_\alpha)$, to $(\alpha_E\oplus 0)\star (0\oplus\bar{\alpha}_E)$. The restriction of the latter map to $qA$, which we denote by $\alpha^\oplus$, is described by $\alpha^\oplus =\alpha_E\pi_A\oplus \bar{\alpha}_E \breve{\pi}_A$. We have, with $\sim$ standing for homotopic,
$$\beta q(\alpha) \sim \beta'q(\alpha) \sim \beta' q(\alpha^\oplus)\sim \beta' q(\alpha_E \pi_A)\oplus \beta' q(\bar{\alpha}_E \breve{\pi}_A)
$$
For the last `$\sim$' we have used the homotopy equivalence $\cK\otimes q(M_2A)\sim \cK\otimes M_2(qA)$. By the uniqueness of the additive homotopy inverse we have that $\beta' q(\bar{\alpha}_E\breve{\pi}_A)\sim \breve{\beta}' q(\bar{\alpha}_E\pi_A)$.
\eproof
We can now apply this procedure to the universal split exact sequence
$$0\to qA\stackrel{j}{\to} QA \to A\to 0$$
and choose for $\alpha, \beta$ the following universal maps $\alpha^u$ and $\beta^u$:
$$\alpha^u =\id_{qA}:qA\to qA,\quad\beta^u =q(\id_{qA})=\id_{q^2A}:q(qA)=q^2 A\to q^2 A.$$
For $\alpha^u_E,\bar{\alpha}^u_E$ we can take $\iota_A,\bar{\iota}_A :A\to QA$. By Proposition \ref{lemext} we can choose $\beta^{u '}: q(QA) \to M_2(q^2A)$ such that $\beta^{u '}q(j)$ is homotopic to $\beta^u =\id_{q^2A}$.\\
From Proposition \ref{propomega} we get a map $\omega^u: qA\to M_2(q^2A)$ such that $\omega^u q(\pi_A)$ is homotopic to $\beta^uq(\alpha^u)=\id_{q^2A}$. In the universal situation here we call this map $\vp_A$. Thus $\vp_A q(\pi_A)$ is homotopic to $\id_{q^2A}$. If we compose $\vp_A$ on the left with $q(\pi_A)$ to $q(\pi_A)\vp_A$, the result is homotopic to $\id_{qA}$. This follows since $q(\pi_A)$ annihilates the terms $\bar{\eta}'\bar{\alpha}_E$ and
$\eta'\bar{\alpha}_E$ in the formula for $\vp_A=\omega^u$ in Proposition        \ref{propomega} and sends $\alpha_E$ to $\id_A$. Thus $\omega^u$ is a homotopy inverse for $q(\pi_A)$. By uniqueness of the homotopy inverse we see that $\omega^u$ must coincide (up to homotopy) with the universal map $\vp_A$ constructed in \cite{CuKK}, \cite{CuGab}.
\bremark
As we see from the discussion above it is more natural to work with $q(\pi_A)$ rather than with $\pi_{qA}$ as in \cite{CuKK}, \cite{CuGab}. This makes part of the arguments quite a bit simpler. Note that we had to argue in \cite{CuKK}, \cite{CuGab} that $q(\pi_A)$ is homotopic to $\pi_{qA}$.\eremark
Consider now again the general situation with homomorphisms $\alpha:qA\to \cK\otimes B$ and $\beta:q(\cK\otimes B)\to C$ as in Proposition \ref{propomega}. We have the homotopy $\beta q(\alpha)\sim \omega q(\pi_{A})$. We have just seen that $q(\pi_A)$ has a homotopy inverse. Thus we see that $\omega$ is uniquely determined up to homotopy by the homotopy classes of $\alpha$ and $\beta$. We can thus make the following definition.
\bdefin\label{defprod}
The product $\alpha\sharp\beta\in KK(A,B)$ of $[\alpha]\in KK(A,B)$ and $[\beta]\in KK(B,C)$ is the homotopy class $[\omega]\in KK(A,C)$.
\edefin
\bremark\label{remhom}
The formula for $\omega$ shows that the product is bilinear with respect to the addition induced by $\oplus$.\\
Any homomorphism $\vp:X\to Y$ induces a quasihomomorphism $(\vp,0)$ and thus an element of $KK(X,Y)$ which we denote by $KK(\vp)$. The composition of $\vp$ with a quasihomomorphism $q(\alpha,\bar{\alpha})$ is $\vp \,q(\alpha,\bar{\alpha})$ or $q(\alpha\vp,\bar{\alpha}\vp)$, respectively. The formula for $\omega$ in \ref{propomega} shows that the $\sharp$-product with $KK(\vp)$ is induced by composition with $\vp$.\\
Since the product will be associative, $KK$ becomes a functor from the category of (separable) C*-algebras to the additive category with objects (separable) C*-algebras and morphism sets $KK(\cdot,\cdot)$.
\eremark
\subsection{Associativity}\label{sua} We follow here the discussion in Section 4 of \cite{CuGen}. Assume that we have elements in $KK(A,B)$, $KK(B,C)$, $KK(C,D)$ represented by homomorphisms $\alpha: qA\to \cK\otimes B$, $\beta: qB \to \cK\otimes C$, $\gamma: qC\to \cK\otimes D$. We use Remark \ref{remstab} to extend $\beta,\gamma$ to homomorphisms from $q(\cK\otimes B), q(\cK\otimes C )$ and we define successively first $E_\alpha\supset B_0$ and $\alpha_E,\bar{\alpha}_E:A\to E_\alpha$ as above, then $\beta':qE_\alpha \to \cK\otimes C$ such that the restriction of $\beta'$ to $qB_0$ is homotopic to $\beta$. We let $C_0$ denote the hereditary subalgebra of $\cK\otimes C$ generated by $\beta'(qE_\alpha)$. Then we define $E_{\beta'}$ as before and get homomorphisms $\beta'_E,\bar{\beta}'_E:E_\alpha\to E_{\beta'}$. We then take $\gamma':qE_{\beta'} \to \cK\otimes D$ such that its restriction to $qC_0$ is homotopic to $\gamma$ and get homomorphisms $\gamma'_E,\bar{\gamma}'_E:E_{\beta'}\to E_{\gamma'}$.

We can now apply Proposition \ref{propomega} to determine the two products $\gamma'\,\sharp\, (\beta' \sharp \alpha)$ and $(\gamma'\sharp\beta' )\, \sharp\, \alpha$. We write here the products in the same order as composition of homomorphisms to make the formulas more intuitive. By the choice of $\beta',\gamma'$ the products $\gamma'\,\sharp\, (\beta' \sharp \alpha)$ and $(\gamma'\sharp\beta' )\, \sharp\, \alpha$ will be homotopic to $\gamma\sharp(\beta\sharp\alpha)$ and $(\gamma\sharp\beta)\sharp\alpha$. By Definition \ref{defprod} the previous products can be described as $\gamma'\sharp \,\omega_1$ and $\omega_2\sharp \,\alpha$ with
\bglnoz
\omega_1
 =q(\beta_E'\alpha_E\oplus\bar{\beta}_E'\bar{\alpha}_E \, ,\, \bar{\beta}_E'\alpha_E\oplus\beta'_E\bar{\alpha}_E )\\
 \omega_2 =q(\gamma_E'\beta'_E\oplus\bar{\gamma}_E'\bar{\beta}'_E \, ,\, \bar{\gamma}_E'\beta'_E\oplus\gamma_E\bar{\beta}'_E )
\eglnoz
We can now apply Proposition \ref{propomega} to both products. By the special form of $\omega_1$, the homomorphisms $\gamma'_E,\bar{\gamma}_E'$ can be composed with the homomomorphisms occuring in the two components of $\omega_1$. Therefore $\gamma'$ extends to $E_{\omega_1}$ and we are in the situation of \ref{propomega}. Second, the two homomorphisms defining $\omega_2$ can be composed with $\alpha_E,\bar{\alpha}_E$ and therefore $\omega_2$ extends to $E_\alpha$. When we apply Proposition \ref{propomega} to $\gamma'\,\sharp\, (\beta' \sharp \alpha)$ and $(\gamma'\sharp\beta' )\, \sharp\, \alpha$ and use the special form of $\omega_1,\omega_2$ we find that in both cases the triple product is given by
$$q(\gamma'_E\beta'_E\alpha_E\oplus\bar{\gamma}'_E\bar{\beta}'_E\alpha_E\oplus
\gamma'_E\bar{\beta}'_E\bar{\alpha}_E\oplus\bar{\gamma}'_E\beta'_E\bar{\alpha}_E \,,\,
\bar{\gamma}'_E\beta'_E\alpha_E\oplus\gamma'_E\bar{\beta}'_E\alpha_E\oplus
\bar{\gamma}'_E\bar{\beta}'_E \bar{\alpha}_E\oplus\gamma'_E\beta'_E\bar{\alpha}_E)
$$
\section{The product in $KK^{op}$}
The formulas in section \ref{secprod} for the definition of the product in $KK$ basically carry over to $KK^{op}$. In fact all homotopies used in section \ref{secprod} for the construction of the product are in fact operator homotopies. As pointed out in Remark \ref{propdeg} this would be enough in order to construct the product for a slightly modified version of $KK^{op}$. However it is not clear if this approach could easily give the result that two $\cK$-stable quasihomomorphisms define the same element of $KK$ iff one can add a degenerate to each of them so that they become operator homotopic (see Corollary \ref{cordeg}).\\
Therefore in this section we will instead work with $KK^{op}$ as defined above and construct the product for $KK^{op}$ again from the start - at the cost of some redundancy. We will use the framework of quasihomomorphisms rather than that of homomorphisms from $qA$. We believe that this also makes the arguments a little more explicit and easier to follow. We use the conventions (such as the notion of a $\cK$-stable quasihomomorphism) of section \ref{secophom} and assume in addition that C*-algebras are $\sigma$-unital whenever necessary. We will assume that all quasihomomorphisms are $\cK$-stable, we will use freely Remark \ref{remstab} and we replace algebras such as $B_0,C_0,D_0$ by $\cK\otimes B,\cK\otimes C,\cK\otimes D$.\\
Here is an adapted version of Proposition \ref{lemext}.
\bprop\label{lemext2}
Let $0 \to \cK\otimes B\stackrel{j}{\to} E \stackrel{p}{\to} A \to 0$ be an exact sequence of C*-algebras with a splitting homomorphism $s:A\to E$ (i.e. $ps=\id_A$) and $(\beta,\bar{\beta}) :\cK\otimes B\to \cK\otimes C$ a $\cK$-stable quasihomomorphism. There is a $\cK$-stable quasihomomorphism $(\beta',\bar{\beta}'):E \to M_2(\cK\otimes C)$ such that the restriction $(\beta'j,\bar{\beta}'j)$ to $\cK\otimes B$ is operator homotopic to $(\beta,\bar{\beta})$.
\eprop
\bproof
Consider the C*-algebra $R$ generated in $M_2(\cM(\cK\otimes C))$ by the matrices in $$\left(\begin{matrix}
 R_1 & R_1R_2\\
 R_2R_1 & R_2
\end{matrix}\right)$$
where $R_1=\beta(\cK\otimes B)$, $R_2=\bar{\beta}(\cK\otimes B)$.\\
Now $\beta,\bar{\beta}$ extend to homomorphisms $\beta^\circ:\cM(\cK\otimes B)\to \cM(R_1)$ and $\bar{\beta}^\circ:\cM(\cK\otimes B)\to \cM(R_2)$. Composing these with the map $E\to \cM(\cK\otimes B)$ we also get homomorphisms
$\beta^E:E\to \cM(R_1)$ and $\bar{\beta}^E:E\to \cM(R_2)$ which, by definition of $R$ we can combine to a homomorphism $\beta^E\oplus \bar{\beta}^E:E\to \cM(R)$ (where $\oplus$ denotes the diagonal sum in $M_2$).\\
Consider also the C*-subalgebra $D$ of $M_2(\cM (\cK\otimes C))$ consisting of matrices of the form
$$\label{D''} d(x)=\left(\begin{matrix}
 \beta^E (x) & 0\\
 0 & \bar{\beta}^E(x)
\end{matrix}\right)\quad x\in E$$
Let $J= R\cap M_2(\cK\otimes C)$. As in the proof of \ref{lemext} $J$ is an ideal in $R$ and also an ideal of $R+\!D$. The quotient $(R+D)/J$ is isomorphic to $M_2(\beta(\cK\otimes B)/\cK\otimes C)+\dot{D}$ where $\dot{D}$ is the quotient of $D$ mod $J$ and consists of the matrices in $\cM(J)/J$ of the form $\left(\scriptsize{\begin{matrix}
 \dot{\beta^E} (x) & 0\\
 0 & \dot{\beta^E} (x)
\end{matrix}}\right)$ with $x\in E$ and where $\dot{\beta^E} (x)$ denotes the image of $\beta^E (x)$ in $\cM (J)/J$ (note that $\beta^E (x)$ and $\bar{\beta}^E (x)$ become equal in the quotient).\\
We can apply Thomsen's extension theorem and lift the multiplier $\left(\scriptsize{\begin{matrix}
 0 & 1\\
 1 & 0
\end{matrix}}\right)$ of $(R+D)/J$ to a selfadjoint multiplier $S$ of $J$ that commutes mod $J$ with $D$. We set $F=e^{\frac{\pi i}{2} S}$ and denote by $\sigma$ the automorphism $\Ad F$ of $\cM (J)$. \\
Consider the pair of homomorphisms $E\to \cM(J)$
\bgl\label{form'2}
\beta'=\left(\begin{matrix}
 \beta^E &0\\
 0 & \bar{\beta}^E(sp)
\end{matrix}\right),\quad \bar{\beta}'=\sigma\,\left(\begin{matrix}
 \beta^E                                                                   \circ(sp)&0\\
 0 & \bar{\beta}^E
\end{matrix}\right)\egl
Note that $\beta'(x)=d(x)+(0\oplus\bar{\beta}^E(sp(x)-x))$. Using the fact that $\sigma$ by definition fixes $d(x)$ mod $J$ and moves $0\oplus \bar{\beta}^E(sp(x)-x))$ to $\beta^E(sp(x)-x) \oplus 0$ mod $J$ we see that $\beta'(x) = \bar{\beta}'(x)$ mod $J$ for $x\in E$.
Therefore the pair $(\beta',\bar{\beta}')$ defines a quasihomomorphism $E\rightrightarrows J$.
Since $spj=0$, the restriction $(\beta'j,\bar{\beta}'j)$
is operator homotopic in $2\times 2$ matrices to $(\beta,\bar{\beta})$.\\
\eproof
As mentioned before, $(\beta',\bar{\beta}')$ is simply the product of $(\beta,\bar{\beta})$ by the projection quasihomomorphism $\kappa_E:E\to \cK\otimes B$.
It is also important to note that in the preceding proof we can choose the same $F$ and $\sigma$ for all splittings $s$. It is also very important that $(\beta',\bar{\beta}')$ is well defined up to operator homotopy, i.e. independent of the choice of $F$ and thus $\sigma$. In fact if $F_1, F_2$ are two choices of the form $e^{ih}$ with $h\in \cK\otimes C$, then $F_1F_2^*$ is in $1+ \cK\otimes C$ and homotopic to $1$.
We now come to the definition of the product in $KK^{op}$ . Let $(\alpha,\bar{\alpha})$ and $(\beta,\bar{\beta})$ be $\cK$-stable quasihomomorphisms from $A$ to $\cK\otimes B$ and from $\cK\otimes B$ to $\cK\otimes C$. Adding a degenerate if necessary we can assume that $\alpha$ and $\bar{\alpha}$ are injective and we can use $C^*(\alpha (A), \bar{\alpha}(A), \cK\otimes B )$  for the $E_\alpha$ of Section \ref{secprod}. For this choice of $E_\alpha$ we again get a split exact sequence $0 \to \cK\otimes B\stackrel{j}{\to} E_\alpha \stackrel{p}{\to} A \to 0$ with two splittings given by $\alpha,\bar{\alpha}$.\\
Let then $(\beta',\bar{\beta}')$ be an extension of $(\beta,\bar{\beta})$ to a $\cK$-stable quasihomomorphism from $E_ \alpha$ as in Proposition \ref{lemext2}.
According to Proposition \ref{propomega} the product of $(\alpha,\bar{\alpha})$ and $(\beta,\bar{\beta})$ must be the quasihomomorphism $A\to \cK\otimes C$ given by the sum of pairs
$\omega = (\beta'\alpha, \bar{\beta}'\alpha)\oplus(\bar{\beta}'\bar{\alpha},\beta'\bar{\alpha})$.\\
Recall that $(\alpha,\bar{\alpha})$ and $(\beta,\bar{\beta})$ represent the same element in $KK^{op}(A,B)$ iff there is a degenerate pair $(\theta,\theta )$ such that $(\alpha\oplus \theta,\bar{\alpha}\oplus \theta)$ is operator homotopic to $(\beta\oplus \theta,\bar{\beta}\oplus \theta)$.
\begin{propdef}\label{propdep}
In the present situation the operator homotopy class of $\omega$ depends only on the operator homotopy classes $[(\alpha,\bar{\alpha})]_{oph}$ and $[(\beta,\bar{\beta})]_{oph}$ of $(\alpha,\bar{\alpha})$ and $(\beta,\bar{\beta})$ of $(\alpha,\bar{\alpha})$ and $(\beta,\bar{\beta})$, respectively. The equivalence class of $\omega$ in $KK^{op}$ also does not change if we add a degenerate to $\alpha$ or $\beta$.\\
We define the product $KK^{op}(A,B)\times KK^{op}(B,C)$ by $[(\alpha,\bar{\alpha})]_{op}\,\sharp\, [(\beta,\bar{\beta})]_{op}= [\omega]_{op}$.
\end{propdef}
\bproof
We have to show that the image in $KK^{op}$ of $\omega = (\beta'\alpha, \bar{\beta}'\alpha)\oplus(\bar{\beta}'\bar{\alpha},\beta'\bar{\alpha})$ does not change if we replace $(\alpha,\bar{\alpha})$ or $(\beta,\bar{\beta})$ by operator homotopic quasihomomorphisms. It is important to note that in the expression for $\omega$, the quasihomomorphism $(\alpha,\bar{\alpha})$ is used only as a pair of single homomorphisms $\alpha$ and $\bar{\alpha}$.

(a) Assume that $(\beta_1,\bar{\beta}_1)$ is elementary operator homotopic to $(\beta_2,\bar{\beta}_2)$ and that $\omega_1,\omega_2$ are the corresponding products. The fact that $(\beta_1,\bar{\beta}_1)$, $(\beta_1,\bar{\beta}_1)$ are elementary operator homotopic means that $\beta_1=\beta_2$ and there is a norm continuous family $[0,1]\ni t\mapsto U_t$ of unitaries in $\cM(\cK\otimes B)$ such that $U_t\bar{\beta}_1(x)-\bar{\beta}_2(x)U_t\in \cK\otimes B$ for all $x,t$ and such that, with $U=U_1$, we have $\bar{\beta}_2=\Ad U \bar{\beta}_1$. If $(\beta_1',\bar{\beta}'_1)$ and $(\beta_2',\bar{\beta}_2')$ are determined by formula \eqref{form'2}, with corresponding inner automorphisms $\sigma_1,\sigma_2$, then $\bar{\beta}_2'= \sigma_2 \Ad (1\oplus U)\sigma_1^{-1}\,\bar{\beta}'_1$. Because $\bar{\beta}_1=\bar{\beta}_2$ mod $\cK\otimes B$, we can actually choose $\sigma_1=\sigma_2$. Since $\sigma_1,\sigma_2$ are of the form $\Ad e^{ih}$ we see that $(\beta_1',\bar{\beta}'_1)$ is elementary operator homotopic to $(\beta_2',\bar{\beta}_2')$. By the definition of $\omega$ this shows that the corresponding products $\omega_1$ and $\omega_2$ are operator homotopic. The compatibility with unitary equivalence is obvious.

(b) Let $t\mapsto U_t$, $t\in [0,1]$ be a norm continuous family of unitaries in $\cM (\cK\otimes B)$ or in $1+\cK\otimes B$ that implements an operator homotopy of $(\alpha,\bar{\alpha})$ with $(\alpha,\Ad U_1\bar{\alpha})$ and write $\bar{\alpha}_{U_t} =\Ad U_t\bar{\alpha}$. Since $\bar{\alpha}_{U_t}(x)$ equals $\bar{\alpha}(x)$ mod $\cK\otimes B$, we have that $\bar{\alpha}_{U_t}(x) \in E_\alpha$ for all $x,t$. Thus, for the definition of the product by $(\beta,\bar{\beta})$, we can take the same $(\beta',\bar{\beta}')$ for all $\bar{\alpha}_{U_t}$. \\
We get that $\beta'\bar{\alpha}_{U_t}=\Ad V_t\beta'\bar{\alpha}$ with $V_t=\beta'^\circ (U_t)$ and
$\bar{\beta}'\bar{\alpha}_{U_t}=\Ad W_t\bar{\beta}'\bar{\alpha}$ with $W_t=\bar{\beta}'^\circ (U_t)$.

We get that,  for $t=1$, the second summand of $\omega$ becomes $(\bar{\beta}'\bar{\alpha}_{U_1},\beta'\bar{\alpha}_{U_1})=(\Ad W_1\bar{\beta}'\bar{\alpha},
\Ad V_1\beta'\bar{\alpha})$. This is unitarily  equivalent to $(\bar{\beta}'\bar{\alpha},
\Ad W_1^*V_1\beta'\bar{\alpha}))$ which in turn is elementary operator homotopic to the second term of the original $\omega$ via the continuous family $t\mapsto W^*_tV_t$.\\
The case of a unitarily equivalent choice for $(\alpha,\bar{\alpha})$ follows similarly.

(c) The formula for $\omega$ shows immediately that the product of $(\alpha,\bar{\alpha})$ or $(\beta,\bar{\beta})$ by a degenerate $(\theta,\theta)$ is homotopic to a degenerate via a rotation in $2\times 2$-matrices. Thus the class in $KK^{op}$ of the product remains unchanged if we add a degenerate to $(\alpha,\bar{\alpha})$ or $(\beta,\bar{\beta})$.
\eproof
The proof of associativity of the product in $KK^{op}$ then follows verbatim the argument in Section \ref{secprod} if we work only with $\cK$-stable quasihomomorphisms, replace the algebras $B_0,C_0,D_0$ by $\cK\otimes B, \cK\otimes C, \cK\otimes D$ and replace homotopy by operator. We therefore have
\bprop
The assignment $([\alpha]_{op}, [\beta]_{op})\mapsto [\alpha]_{op}\sharp [\beta]_{op}$ defines an associative bilinear product $KK^{op}(A,B)\times  KK^{op}(B,C)\to KK^{op}(A,C)$.
\eprop
\bremark\label{remstable}
The functor $KK^{op}$ is stable in the sense that the natural inclusion map $j_A: A\to \cK\otimes A$ defines an invertible element in $KK^{op}(A,\cK\otimes A)$ for each $A$. This follows from the fact that $\id_\cK\otimes j_A$ is unitarily equivalent to $\id_{\cK\otimes A}$.
\eremark

\section{Split exactness}
A functor $F$ from the category of C*-algebras to an additive category $C$ is called split exact if, for every exact sequence $0\to J\stackrel{j}{\to} E\stackrel{p}{\to} A\to 0$ of C*-algebras with splitting $s:A\to E$, the induced map $F(j)+F(s)$ gives an isomorphism $F(J)\oplus F(A) \to F(E)$ in $C$  (the inverse isomorphism is then necessarily given by $(F(\id)-F(sp))\oplus F(p)$.
\blemma\label{lemspli}
Let $0\to J\stackrel{j}{\to} E\stackrel{p}{\to} A\to 0$  be an exact sequence of C*-algebras with splitting $s:A\to E$. Consider the quasihomomorphism $\kappa_E =(\id_E, sp)$ from $E$ to $J$. For the compositions of this quasihomomorphism with the homomorphisms $s$ and $j$ we have the following operator homotopies:
$$
\kappa_E\,s \sim d_1 \quad \kappa_E \,j \sim \id_J \quad j\,\kappa_E \oplus sp \,\sim\, \id_E\oplus d_2
$$
where $d_1,d_2$ are the degenerates $d_1=(s,s)$ and $d_2=(sp,sp)$.
\elemma
\bproof
The first two homotopies are obvious and in fact even identities. The quasihomomorphism $\kappa_E \oplus sp$ is given by the pair
$$
\left(\left( \begin{matrix}
 \id_E &0\\
 0 & sp
\end{matrix} \right)\, ,
\left( \begin{matrix}
 sp &0\\
 0 & 0
\end{matrix} \right)\right)
$$
This is operator homotopic to $\id_E\oplus (sp,sp)$ via a rotation of the term $sp$ in the second component to the lower right corner.
\eproof
\bprop\label{propsex}(cf. \cite[Proposition 2.1]{CuKK})
The functors $KK^{op}$ and $KK$ are split exact.
\eprop
\bproof
Let as before $0\to J\stackrel{j}{\to} E\stackrel{p}{\to} A\to 0$ be an exact sequence of C*-algebras with splitting $s:A\to E$. Let$KX$ be one of the functors $KK$ or $KK^{op}$. We construct elements in $KX(E,J\oplus A)$ and in $\ KX(J\oplus A,E)$ which are inverse to each other in $KX$.\\
We need some notation. We denote by $pr_A, pr_J$ the two projections from $A\oplus J$ to $A,J$ and by $j_A, j_J$ the inclusions of $A,J$ into the first and second components of $A\oplus J$. With this notation we define the quasihomomorphism $\vp$ from $E$ to $A\oplus J$ by $\vp = j_A\,p + j_J\,\kappa_E$ and the homomorphism $\psi: A\oplus J \to M_2(E)$ by $\psi = s\, pr_A \oplus j\, pr_J$ where `$\oplus$' denotes the direct sum in $M_2(E)$.\\
Since $\psi$ is a homomorphism, the $KX$-products are again just the composition with $\psi$. We denote the products by $\psi\vp$ and $\vp\psi$ and use the symbol $\sim$ for equivalence in $KK^{op}$ (i.e. addition of degenerates + operator homotopy ). For $\psi\vp$ we get using Lemma \ref{lemspli} that $\psi\vp=sp\oplus j\kappa_E \sim \id_E$. For $\vp\psi$ we get the matrix
$$
\left( \begin{matrix}
 j_A\,ps\,pr_A + j_J\kappa_E s\,pr_A &0\\
 0 & j_A pj \,pr_J + j_J \kappa_E j\, pr_J
\end{matrix} \right)
$$
Since using Lemma \ref{lemspli} again $ps=\id_A,\, \kappa_E s\,\sim 0,\, pj=0,\, \kappa_E j\sim \id_J$ we get that $\psi\vp$ is equivalent in $KX$ to $\id_{A\oplus J}$.
\eproof
\section{Homotopy invariance of $KK^{op}$ and a criterion for the equivalence of two quasihomomorphisms in $KK$}
Given a C*-algebra $A$ we denote by $CA$ the algebra $A[0,\infty)$ of continuous $A$-valued functions on the interval $[0,\infty)$ that vanish at infinity. We denote by $\ev_t: CA \to A$ the evaluation homomorphism that evaluates at $t\in [0,\infty)$. The proof of the following proposition is inspired by but, at least for $KK^{op}$, shorter than Ralf Meyer's proof of the fact that any split exact and stable functor is homotopy invariant (on many natural categories of algebras) \cite[Lemma 3.26]{CMR}.
\bprop\label{propev0}
In $KK^{op}(CA,A)$ we have $KK^{op}(\ev_0)=0$.
\eprop
\bproof
Denote by $p_k$ the projection onto the one-dimensional subspace of $\ell^2(\Zz)$ spanned by the $k$-th element in the standard basis.
Given $n\in \Nz$, consider the homomorphisms $\vp_n,\bar{\vp}_n:CA\to \cK(\ell^2\Zz)\otimes A$ given by
$$\vp_n(f) =\sum_{k=0}^\infty f(k2^{-n})p_k\qquad \bar{\vp}_n (f) =\sum_{k=0}^\infty f((k+1)2^{-n})p_k.$$
Then $\vp_n(f),\bar{\vp}_n(f)$ are both compact, i.e. in $\cK(\ell^2\Zz)\otimes A$, and the pair $(\vp_n,\bar{\vp}_n)$ defines a quasihomomorphism $CA\to \cK(\ell^2\Zz)\otimes A$. The shift operator $U$ on $\ell^2(\Zz)$ transports $\bar{\vp}_n$ to $\vp^0_n$ defined by $$\vp^0_n(f)=\sum_{k=0}^\infty f((k+1)2^{-n})p_{k+1}.$$ Note that unlike in the definition of $\bar{\vp}_n$ we use here $p_{k+1}$ rather than $p_k$.
Replacing $U$ by the $2\times 2$-matrix $U\oplus U^*$ we can connect $U$ continuously to 1 and get a homotopy between $\bar{\vp}_n$ and $\vp^0_n$. Since $\bar{\vp}_n(f)$ and $\vp^0_n(f)$ are already compact for each $f$ this gives an operator homotopy between $(\vp_n,\bar{\vp}_n)$ and $(\vp_n,\vp^0_n)$.
Obviously $\vp_n =\ev_0\oplus \vp_n^0$ so that $KK^{op}(\vp_n,\bar{\vp}_n)=K^{op}(\vp_n,\vp_n^0)= KK^{op}(\ev_0,0)$ (we denote here by $KK^{op}(\alpha,\bar{\alpha})$ the element of $KK^{op}$ defined by a quasihomomorphism $(\alpha,\bar{\alpha})$).\\
We have that $\vp_{n+1} = \vp_n'\oplus \psi_n$ and $\bar{\vp}_{n+1} = \bar{\vp}_n'\oplus \psi_n$ where
\bglnoz\vp'_n(f) =\sum_{k=0}^\infty f(k2^{-n})p_{2k}\qquad \bar{\vp}'_n (f) =\sum_{k=1}^\infty f((k+1)2^{-n})p_{2k}\\ \psi_n(f) = \sum_{k=0}^\infty f((2k+1)2^{-(n+1)})p_{2k+1}\eglnoz
(note that for $\vp_n',\bar{\vp}_n'$ we use $p_{2k}$ rather than $p_k$) and $\vp'_n,\bar{\vp}_n'$ are obviously unitarily equivalent to $\vp_n,\bar{\vp}_n$. Thus $KK^{op}((\vp_n,\bar{\vp}_n))=KK^{op}((\vp'_n,\bar{\vp}'_n)\oplus (\psi_n,\psi_n))=KK^{op}(\vp_{n+1},\bar{\vp}_{n+1})$.\\
Now every $f\in CA$ is uniformly continuous on $[0,\infty)$ and therefore $\|\vp_n(f)-\bar{\vp}_n(f)\|\to 0 $ for $n\to \infty$. For each $j\in \Nz$ we can therefore define a quasihomomorphism $\Phi_j:CA\to \cK(\ell^2(\Zz\times \Zz)\otimes A$ (where the first copy of $\Zz$ stands for the index $k$ and the second one for the index $n$) by
$$\Phi_{j} = \left(\sum_{n\geq j} \vp_n\otimes e_n,\,\sum_{n\geq j}\bar{\vp}_n\otimes e_n\right)$$
where $e_n$ is the projection onto the $n$-th basis vector in the second copy of $\ell^2\Zz)$.
Writing $\Psi_j$ for the degenerate quasihomomorphism $\left(\sum_{n\geq j} \psi_n\otimes e_n,\sum_{n\geq j} \psi_n\otimes e_n\right)$ we get that $\Phi_j  \oplus \Psi_j$ is unitarily equivalent to $\Phi_{j+1}$ and thus $\Phi_j$ is equivalent in $KK^{op}$ to $\Phi_{j+1}$. Summarizing we see for $j=0$ that $KK^{op}(\vp_0,\bar{\vp}_0)+KK^{op}(\Phi_0)=KK^{op}(\vp_0,\bar{\vp}_0)+KK^{op}(\Phi_1)=KK^{op}(\Phi_0)$ which implies that $KK^{op}(\ev_0)=KK^{op}(\vp_0,\bar{\vp}_0)=0$.
\eproof
\btheo\label{main}
There is a natural map $KK^{op}(A,B)\to KK(A,B)$ which is an isomorphism for each $A$ and $B$ - assuming that $A$ is separable. In particular $KK^{op}$ is homotopy invariant in both variables.
\etheo
\bproof
Every $\cK$-stable quasihomomorphism $(\vp,\bar{\vp}):A\to \cK\otimes B$ defines a homomorphism $q(\vp,\bar{\vp}):qA\to \cK\otimes B$. If two quasihomomorphisms are operator homotopic then the corresponding homomorphisms from $qA$ are homotopic. Also any degenerate quasihomomorphism $(\theta,\theta))$ is sent to 0 in this correspondence. Since by Kasparov's stabilization theorem every homomorphism $qA\to \cK\otimes B$ can be realized by a $\cK$-stable quasihomomorphism $A\rightrightarrows \cM (\cK\otimes B) \rhd \cK\otimes B$, the resulting map $KK^{op}(A,B)\to KK(A,B)$ is well defined and surjective.\\
A $\cK$-stable quasihomomorphisms $(\alpha,\bar{\alpha})$ from $A$ to $\cK\otimes B$ represents 0 in $KK(A,B)$ iff there is a quasihomomorphism (which we may assume to be  $\cK$-stable) $(\vp,\bar{\vp}):A \to \cK\otimes CB$ such that $(\vp,\bar{\vp})\sharp \ev_0=(\alpha,\bar{\alpha})$. But then $KK^{op}(\alpha,\bar{\alpha})=KK^{op}(\vp,\bar{\vp})KK^{op}(\ev_0)=0$ by Proposition \ref{propev0}. This shows that the map $KK^{op}(A,B)\to KK(A,B)$ is injective.
\eproof
This theorem has been established by Kasparov \cite{KasExt}. An elegant proof is also due to Skandalis \cite{SkanOp}.
We obtain the following important consequence  of Theorem \ref{main}
\bcor\label{cordeg}
Two $\cK$-stable quasihomomorphisms $(\alpha,\bar{\alpha})$ and $(\beta,\bar{\beta})$ from $A$ to $\cK\otimes B$ represent the same element in $KK(A,B)$ iff there is a degenerate $(\theta,\theta)$ such that $(\alpha,\bar{\alpha})\oplus (\theta,\theta)$ is operator homotopic to $(\beta,\bar{\beta})\oplus (\theta,\theta)$.\\
\ecor
Assume now that $F$ is any functor from the category of separable C*-algebras to an additive category which is split exact and stable in the sense that $F$ transforms the natural embedding of $A$ into $\cK\otimes A$ into an isomorphism for each $A$. We know that $KK^{op}$ is split exact (Proposition \ref{propsex}) and stable (Remark \ref{remstable}). Stability (in fact already $M_2$-stability) implies that $F$ is invariant under unitary equivalence and split exactness implies that $F$ is additive. In the paragraph before Proposition \ref{propomega} we have seen that we can associate, with any $\cK$-stable quasihomomorphism $(\alpha,\bar{\alpha})$ between C*-algebras $A$ and $B$, a split exact sequence $0\to\cK\otimes B \to E_\alpha \to A\to 0$ and a quasihomomorphism $(\alpha_E,\bar{\alpha}_E):A\rightrightarrows E_\alpha \triangleright \cK\otimes B$ such that $q(\alpha_E,\bar{\alpha}_E)=q(\alpha,\bar{\alpha})$.  Then by split exactness and stability $F(\alpha_E)-F(\bar{\alpha}_E)$ sends $F(A)$ to $F(\cK\otimes B)=F(B)$. We denote this morphism by $F(\alpha,\bar{\alpha}):F(A)\to F(B)$. See \cite[3.1.1]{CMR} for more details on this construction. Obviously $F(\theta,\theta)=0$ for a degenerate quasihomomorphism $(\theta,\theta)$. We can then apply the argument in the proof of Proposition \ref{propev0} to $F$ in place of $KK^{op}$ and get that $F(\ev_0)=0$.
\bprop\label{prophominv} Let $F$ be as above and
let $A$ be a separable C*-algebra. Let $A[0,1]$ be the C*-algebra of continuous $A$-valued functions on the interval $[0,1]$ and $\ev_1,\ev_0$ the homomorphisms $A[0,1]\to A$ given by evaluation at 1 and 0. The kernel of $\ev_1$ is isomorphic to the algebra $CA$ considered above and $\ev_0$ corresponds, under this isomorphism, to the homomorphism $\ev_0$ used above.\\
We have that $F(\ev_0)=F(\ev_1)$, i.e. $F$ is homotopy invariant.
\eprop
\bproof
We have the split exact sequence $0\to CA\stackrel{j}{\to} A[0,1]\stackrel{\ev_1}{\to} A\to 0$ with splitting $s:A\to A[0,1]$ given by $s(x)=x1$. By split exactness we have $F(A[0,1])\cong F(CA)\oplus F(A)$ and $F(\ev_0)=F(\ev_0|_{CA}) + F(\ev_1)=F(\ev_1)$ since $F(\ev_0|_{CA})=0$.
\eproof

\section{Stable uniqueness}\label{secstab}
In this section we derive the stable uniqueness theorem of Dadarlat-Eilers \cite{DadEil} from Corollary \ref{cordeg}. This theorem was originally derived in \cite{DadEil} from a version of that corollary, too. A somewhat more direct way to do this - even in the equivariant case - was found in \cite{GabSzab}. Using a key idea from \cite[Lemma 4.3]{GabSzab} but arranging the rest of the argument in a slightly different way, we get below a short proof of the stable uniqueness theorem (in the non-equivariant case).
\bdefin
Let $\vp,\psi: A\to \cM(B)$ be homomorphisms. We say that $\vp$ is asymptotically $B$-inner equivalent to $\psi$, if there is a continuous path $[0,\infty)\ni t\mapsto U_t$ where the $U_t$ are unitaries in $1+B$ such that such that $U_t\vp(a)U_t^*-\psi (a)\in B$ for all $t$ and such that $U_t\vp(a)U_t^*$ tends to $\psi (a)$ as $t\to \infty$ for each $a\in A$.
\edefin
Note: what we call asymptotically $B$-inner equivalent is called `properly asymptotically unitarily equivalent' in \cite{DadEil}.\\
Let $t$ be a parameter in $\Nz$ or in $[0,\infty)$ and let $(a_t)$ be a family of elements in a C*-algebra $A$. We will consider the exponential function $\exp (x)$ represented by the power series $\sum (1/n!)x^n$. Given $a\in A$ the commutator $\delta_a =[a,\,\cdot]$ satisfies the derivation rule $\delta_a(x^n)=\sum_{k=0}^{n-1} x^k\delta_a(x)x^{(n-1)-k}$ whence $\|\delta_a(x^n)\|  \leq n\|\delta_a(x)\|\|x\|^{n-1}$. If we apply this latter inequality we get the well known fact that the difference $\exp (x+a_t)-\exp (x)\exp(a_t)$ tends to 0, if the commutator $[a_t,x]$ tends to 0 for $t\to \infty$.
\blemma\label{lem9.2}
Let $\vp :A \to \cM (B)$ be a homomorphism from the separable C*-algebra $A$ to the multipliers of the $\sigma$-unital C*-algebra $B$. Let $x$ with $x^*=-x$ in $\cM(B)$ such that $[x,\vp (A)]\subset B$. Let $U=\exp x$ and $(h_t)$ a commutative approximate unit for $B$ which is quasicentral for $\vp (A)$ and for $C^*(x)$. Consider the unitaries $V_t=\exp h_txh_t$ in $1 +B$. Then $V_t^*U - \exp (x-h_txh_t) $ tends to 0 and $[V_t^*U,\vp (a)]$ tends to 0 for each $a\in A$.
\elemma
\bproof
By assumption the commutator $\vp(a)x-x\vp(a)$ is in $B$. Thus we get
$$\vp(a)x-x\vp(a)\sim (\vp(a)x-x\vp(a))h_t^2 \sim \vp(a)xh_t^2-xh_t^2\vp(a)\sim [\vp(a),h_t x h_t]$$
where $\sim$ means that the difference of the two expressions tends to 0 for $t\to \infty$. Taking the difference we get that $[\vp(a),x-h_txh_t]\sim 0$.\\
Applying the remark on the exponential series above we get that $V_t^*U\sim \exp(x-h_txh_t)$. Since we have seen before that $[x-h_txh_t,\vp(a)]\sim 0$ for each $a\in A$ we get that $[V_t^*U,\vp(a)]\sim 0$.
\eproof
\blemma\label{lem9.3}
Let $(\vp,\bar{\vp}):\, A\rightrightarrows \cM(B)\rhd B$ be a quasihomomorphism. Assume that $U=\exp(x)$ with $x^*=-x\in \cM(B)$ induces an elementary operator homotopy between $(\vp,\bar{\vp})$ and $(\vp,\vp)$, i.e. that $\Ad U\vp = \bar{\vp}$. Then there is a continuous family $[0,\infty ) \ni t\mapsto U_t$ of unitaries in $1+B$ such that $\Ad U_t \vp(a)$ tends to $\bar{\vp}(a)$ as $t\to \infty$ for each $a\in A$.
\elemma
\bproof By Lemma \ref{lem9.2} there is such a family given by $\exp (h_t x h_t)$ for $t\in \Nz$.  We get a continuous family for $t\in [0,\infty)$ by interpolating linearly between $h_t x h_t$ and $h_{t+1} x h_{t+1}$.\eproof
\btheo
(a) A $\cK$-stable quasihomomorphism $(\vp,\bar{\vp}):\, A\rightrightarrows \cM(\cK\otimes B)\rhd \cK\otimes B$ represents 0 in $KK(A,B)$ if and only if there is a homomorphism $\theta:A\to \cM(\cK\otimes B)$ such that such that $\bar{\vp}\oplus \theta$ is asymptotically $\cK\otimes B$-inner equivalent to $\vp\oplus \theta$.\\
(b) Two $\cK$-stable quasihomomorphisms $(\alpha,\bar{\alpha}):\, A\rightrightarrows \cM(\cK\otimes B)\rhd \cK\otimes B$ and $(\beta,\bar{\beta}):\, A\rightrightarrows \cM(\cK\otimes B)\rhd \cK\otimes B$ represent the same element in $KK(A,B)$ iff there is a homomorphism $\theta:A\to \cM(\cK\otimes B)$ and a unitary $W\in \cM(\cK\otimes B)$ such that $(\alpha\oplus \theta,\bar{\alpha} \oplus\theta)$ is asymptotically $\cK\otimes B$-inner equivalent to $(\Ad W(\beta\oplus \theta),\Ad W(\bar{\beta} \oplus\theta))$.
\etheo
\bproof
(a) If $[0,\infty)\ni t\mapsto U_t\in 1+\cK\otimes B$ is a unitary path such that $U_0=1$ and $\Ad U_t(\alpha) \to \bar{\alpha}$, then $[0,1]\ni s\mapsto \Ad U_{1/s}(\alpha)$ gives a homotopy between $\alpha$ and $\bar{\alpha}$, cf.\cite[Lemma 3.1]{DadEil}.\\
Conversely, if $KK(\alpha,\bar{\alpha})$ is 0, then there is $\theta$ and an elementary operator homotopy between $(\alpha\oplus\theta,\bar{\alpha}\oplus\theta)$ and  $(\alpha\oplus\theta,\alpha\oplus\theta)$, i.e. there is a continuous family $[0,1]\ni t\mapsto U_t$ of unitaries in $\cM(\cK\otimes B)$ that commute mod $\cK\otimes B$ with $(\alpha\oplus\theta) (A)$ and are such that $U_0=1$ and $\Ad U_1(\alpha\oplus\theta) =\bar{\alpha}\oplus\theta$. But then $U_1$ is a finite product of unitaries $\exp x_i$ where $x_i^*=-x_i$ and the $x_i$ commute mod $\cK\otimes B$ with $(\vp\oplus\theta)(A)$. We can apply Lemma \ref{lem9.3} to each of the $\exp(x_i)$ and then take the product of the resulting asymptotic unitary paths in $1+\cK\otimes B$.\\
(b) The case of the equivalence of two quasihomomorphisms follows similarly applying Lemma \ref{lemuni}.
\eproof

\end{document}